\documentclass[a4paper,12pt]{amsart}
\usepackage{amssymb,amsmath,amsthm}
\usepackage{graphicx}

\usepackage{amsfonts}
\usepackage{latexsym}
\usepackage{color}
\usepackage[english]{babel}
\usepackage{hyperref}
\usepackage{hyperref}
\usepackage{verbatim} 
\usepackage{bm} 
\usepackage{xcolor}
\usepackage{graphicx} 
\usepackage{subcaption}
\usepackage[font=footnotesize,labelfont=bf]{caption}
\usepackage{mathrsfs}
\definecolor{darkgr}{rgb}{0.0, 0.62, 0.42}

\numberwithin{equation}{section}

\newtheorem{lemma}{Lemma}[section]
\newtheorem{theorem}[lemma]{Theorem}

\newtheorem{remark}[lemma]{Remark}

\newtheorem{definition}[lemma]{Definition}

\def\mass{\mathcal M \kern -2pt a \kern-2pt d}

\def\im{{\rm i}}

\def\sleq{\preceq}
\def\erremu{3R/4}
\def\erremuu{R/2}

\def\cC{{\mathcal{C}}}

\def\cG{{\mathcal{G}}}

\def\cO{{\mathcal{O}}}

\def\cU{{\mathcal{U}}}
\def\cR{{\mathcal{R}}}

\def\cP{{\mathcal{P}}}
\def\cT{{\mathcal{T}}}

\newcommand{\R}{\mathbb R}
\newcommand{\C}{\mathbb C}
\newcommand{\Z}{\mathbb Z}

\newcommand{\T}{\mathbb T}

\def\ac{{\mathcal A}\kern-.7pt\ell\kern-.9pt\mathcal{S}}

\begin{document}

\title[A Nekhoroshev theorem for perturbations of BO]{A Nekhoroshev theorem for some perturbations of the 
  Benjamin-Ono equation with initial data close to finite gap tori}
\date{December 5, 2023}

\author[D.Bambusi]{Dario Bambusi}
\address{Dipartimento di Matematica, Universit\`a degli Studi di Milano, Via Saldini 50, I-20133
Milano} \email{ \tt{dario.bambusi@unimi.it}}
\author[P.G\'erard]{Patrick G\'erard}
\address{Universit\'e Paris--Saclay, Laboratoire de Math\'ematiques d'Orsay,  CNRS,  UMR 8628,  F-91405 Orsay} \email{\tt {patrick.gerard@universite-paris-saclay.fr}}


\maketitle

\begin{abstract} We consider a perturbation of the Benjamin Ono
  equation with periodic boundary conditions on a segment.  We consider
  the case where the perturbation is Hamiltonian and the corresponding
  Hamiltonian vector field is analytic as a map form energy
  space to itself. Let $\epsilon$ be the size of the
  perturbation. We prove that for initial data close in energy norm to
  an $N$-gap state of the unperturbed equation all the actions of
  the Benjamin Ono equation remain $\cO(\epsilon^{\frac{1}{2(N+1)}})$
  close to their
  initial value for times exponentially long with  $\epsilon^{-\frac{1}{2(N+1)}}$.
\end{abstract} \noindent


\medskip

\noindent



\section{Introduction}\label{intro}

In this paper we prove a Nekhoroshev type theorem for
perturbations of the Benjamin-Ono (BO) equation on the torus. It is
well known that the BO equation is a Hamiltonian integrable system,
but only recently it has been understood how to introduce Birkhoff
coordinates \cite{GK21}, namely some regularized action angle variables, for such
a system with periodic boundary conditions. We will call finite gap states, functions $u$ which, when
written in Birkhoff coordinates, have only a finite number of
coordinates different from zero --- see  \cite{GK21} for a justification of this denomination.  Here we consider a small Hamiltonian
perturbation of the BO equation and prove that, for initial data close
to finite gap initial states, all the actions of the
BO equation remain close to their initial value for times which are
exponentially long with a power of $\epsilon^{-1}$, $\epsilon$ being
the size of the perturbation. The main point is that we are not
confined neither to small amplitude initial data nor to neighbourhood
of nonresonant tori: as far as we know this is
the first result of this kind for PDEs. Our result however has
a strong limitation: we deal only with perturbations whose
Hamiltonian vector field is an analytic map from the phase space to itself, and
in view of the fact that the Poisson structure of the BO equation is
the Gardner's one, namely $\partial_x$, this means that the
Hamiltonian must depend on some antiderivative of the state $u$ and
typically the vector field is non local. The main example of
Hamiltonian perturbation that we can treat is given by $\epsilon P(u)$, with
\begin{equation}
  \label{pertu}
P(u):= \frac{1}{2\pi}\int_0^{2\pi} F(x,\partial_x^{-1}u)dx \ ,
\end{equation}
where $F(x,y)$ is analytic in the second variable and continuous in
the first one.  We also cover Hamiltonian versions  of
 perturbations
dealt with by \cite{Gas22}. We will comment further
on this point after the statement of the main theorem.

Notwithstanding this defect, we think that the fact that exponentially
long stability estimates hold for some perturbations of the BO
equation, close to finite dimensional tori could spread some light on
the behaviour of perturbed integrable Hamiltonian systems. In
particular this has to be confronted with the paper \cite{BouKal} in
which a mechanism of possible instability of infinite dimensional tori
was exhibited.

Finally we recall that all known Nekhoroshev type results for PDEs
can be split into two types: results dealing with neighbourhoods of
nonresonant tori (see \cite{Bou00, BFG20, BG21, KM21}) --- which are
essentially stability results for such tori ---, and results dealing with
small amplitude solutions, mainly in perturbations of nonresonant
systems (see e.g. \cite{Bou96,Bam03,BG06,BDGS,BFM22} and references
therein). We also mention the paper
\cite{Bam99} in which a Nekhoroshev type result controlling also
neighbourhood of resonant tori has been proved, for small amplitude solutions in
perturbations of the integrable NLS.

The proof of the main result of the present paper is a variant of the
proof of \cite{Bam99} (see also \cite{BG93}) which is based on
Lochak's proof of Nekhoroshev's theorem \cite{Loc92,Loc2}. We emphasize
that this is made possible by the fact that, as shown by the works
\cite{GK21, GKT23}, the Hamiltonian of the BO equation, written in
Birkhoff coordinates, has a very simple and explicit structure and this
allows us to verify all the geometric properties needed in order to apply
Lochak's technique.

{\it Acknowledgements} The present research was founded by the PRIN project
  2020XB3EFL Hamiltonian and dispersive PDEs. It was also supported by
GNFM. The second author is supported by the project ANR-23-CE40-0015-01.

\section{Main result}\label{main}

Consider the Benjamin-Ono (BO) equation
\begin{equation}
  \label{BO}
\partial_tu={\rm H}\partial_x^2 u-\partial_x(u^2)\ ,\quad
x\in\T:=\R/2\pi\Z\ ,\quad t\in\R 
\end{equation}
where $u$ is real valued and ${\rm H}$ denotes the Hilbert transform defined
by
$$
{\rm H}u(x):=-i\sum_{n\not=0} {\rm sgn}(n)\hat u_n{\rm e}^{inx}\ ,
$$ 
here $\hat u_n:=\frac{1}{2\pi}\int_0^{2\pi}u(x){\rm e}^{-\im nx }dx$ are the
Fourier coefficients of $u$. 
Equation \eqref{BO} is Hamiltonian with the Hamiltonian
function
\begin{equation}
  \label{Hu}
H_{BO}(u):=\frac{1}{2\pi}\int_0^{2\pi}\left(\frac{1}{2}(\left|\partial_x\right|^{1/2}u)^2-\frac{1}{3}u^3\right)dx\ ,
\end{equation}
where $\left|\partial_x\right|$ is the Fourier multiplier by
$\left|n\right|$ and the Poisson tensor is Gardner's one, namely the Hamiltonian equation associated to a Hamiltonian function $H$ is
\begin{equation}
  \label{Hamilton}
\dot u=\partial_x\nabla H(u)\ ,
\end{equation}
where $\nabla$ denotes the $L^2$ gradient defined by $\langle \nabla
H(u),h \rangle=dH(u)h$, for every $ h\in C^{\infty}(\T) $. It is well
known that Equation \eqref{BO} is integrable. The Birkhoff
coordinates have been introduced in \cite{GK21}, and further studied
in \cite{GKT20, GKT22, GKT23, GKTana, GKT21}. To recall the result of these papers consider the space
$H^s_{r,0}\subset H^s(\T,\R)$ of the real valued functions  $u$ of
class $H^s$ having zero mean value and the space $h^{s}_+$ of the
sequences $(\xi_n)_{n\geq 1}$, $\xi_n\in\C$, such that
$$
\left\|{\xi}\right\|_s^{2}:=\sum_{n\geq 1}n^{2s}\left|\xi_n\right|^2<\infty\ .
$$
We also need the notation $H^s_0(\T,\C)$ for $H^s$ functions with zero mean value.
\begin{theorem}
[G\'erard-Kappeler-Topalov, \cite{GKT23, GKTana, GKT21}] 
  There exists a map
  \begin{align*}
    \Phi: \bigsqcup_{s>-1/2}H^s_{r,0}&\to\bigsqcup_{s>-1/2}h_+^{s+1/2},\\
    u&\mapsto
    \xi(u):=(\xi_n(u))_{n\geq 1}
  \end{align*}
so that the following properties hold for any $s>-1/2$.
\begin{itemize}
\item[1]. $\Phi:H^s_{r,0}\to h_+^{s+1/2} $ is a diffeomorphism;
  $\Phi$ and its inverse $\Phi^{-1}$, map bounded sets to bounded sets.

  \item[2]. The map is symplectic, in the sense that in terms of the
    variables $\xi_n$ the Hamilton equations \eqref{Hamilton} take the
    form
    \begin{equation}
      \label{symle}
\dot\xi_n=\im\frac{\partial H}{\partial\bar \xi_n}\ ,\quad 
\dot{\overline{ \xi_n}}=-\im\frac{\partial H}{\partial \xi_n}\ .
    \end{equation}
        \item[3]. In terms of the variables $\xi$ the BO Hamiltonian
          takes the form
          \begin{align}
                   \label{BOH}
H_{BO}(u(\xi))&= H_2-H_4
\\
\label{BOH1}
H_2&=\sum_{n\geq 1}n^2|\xi_n|^2\ ,\quad H_4=\sum_{n\geq
  1}[s_n(|\xi|^2)]^2\ ,
\\
&s_n(|\xi|^2):=\sum_{k\geq n}|\xi_k|^2\ .
            \end{align}
            \end{itemize}
\end{theorem}

In \cite{GKTana, GKT21} it was also proved that the map $\Phi$ is  real analytic, in the sense that, for any $s>-1/2$,  $\Phi $ extends
to a holomorphic map from a complex neighbourhood of $H^s_{r,0}$ in $H^s_0(\T,\C)$ to $h^{s+1/2}:=h^{s+1/2}_+\oplus h^{s+1/2}_+$. A similar result holds for $\Phi^{-1}$. This is crucial for our application.
\vskip 0.25cm
{\it From now on we will work in the energy space $H^{1/2}_{r,0}$.}
\vskip 0.25cm
To state our main result, for any integer $N\geq 1$, consider the set
of $N$ gap states, namely
\begin{equation}
  \label{fg}
\cG^N:=\left\{u\in H^{1/2}_{r,0}\ : |\xi_n(u)|\ne
  0,\  1\leq n\leq N\ ,\ |\xi_n(u)|=0,\ 
n>N\right\}\ .
  \end{equation}
  The positive quantities $\gamma_n (u):=|\xi_n(u)|^2, n=1,\dots, N$, are called the gaps of the state $u\in \cG^N$. \\
In section 7, remark 7.2(i) of \cite{GK21}, it is proved that $\cG^N$ is a dense open subset of the symplectic $2N$-dimensional manifold
$$\cU_N:=\left \{ \sum_{j=1}^N (P_{r_j}(x+\alpha _j)-1)\ ,\ (r_1,\dots ,r_N)\in (0,\infty )^N\ ,\ (\alpha_1,\dots ,\alpha_N)\in \T^N\right \} \ ,$$
where $P_r$ denotes the usual Poisson kernel
$$P_r(x):=\frac{1-r^2}{1-2r\cos x+r^2}\ .$$
\begin{theorem}
  \label{main1}
Consider a Hamiltonian system with Hamiltonian function
\begin{equation}
  \label{hperturbata}
H:=H_{BO}+\epsilon P\ ,
\end{equation}
 and assume that there exists $N\geq 1$ such that  $P$
 extends to a real analytic function on a neighbourhood of $\cG^N$,
 which is bounded on
 bounded sets. Also assume  that its Hamiltonian vector field $X_P$
 extends to a real analytic function from an open neighbourhood of $\cG^N\subset H^{1/2}_{r,0}$
 to $H^{1/2}_{r,0}$ 
and it is is bounded on bounded sets. Fix two positive parameters $0<E_m<E_M$,
 then the following holds true: there exist constants
 $\epsilon_*,C_1,...,C_5$, independent of $\epsilon$ s.t. if
 $|\epsilon|<\epsilon_*$ and the initial datum $u_0$ fulfills
\begin{align}
  \label{ini.1}
  E_m\leq |\xi_n(u_0)|^2\leq E_M\ ,\quad \forall n\leq N\ ,
  \\
  \label{ini.2}
  \sum_{n\geq N+1}n^2|\xi_n(u_0)|^2<C_1\epsilon^{1/2(N+1)}\ ,
\end{align}
then along the flow of the Hamiltonian system \eqref{hperturbata} one
has
\begin{align}
  \label{stime.1}
  \left||\xi_n(t)|^2-|\xi_n(0)|^2\right|\leq C_2\epsilon^{1/2(N+1)}\ ,
  \\
  \sum_{n\geq N+1}n^2|\xi_n(t)|^2\leq C_3 \epsilon^{1/2(N+1)}
  \end{align}
for all times $t$ fulfilling
\begin{equation}
  \label{stima.exp}
|t|\leq C_4 \exp\left(\frac{C_5}{\epsilon^{1/2(N+1)}}\right)\ .
  \end{equation}
  \end{theorem}
A family of examples of perturbations $P$ fulfilling our assumptions
is 
\begin{equation}
  \label{exe.1}
P(u):=\int_{0}^{2\pi}F(x,\partial^{-1}_xu(x))dx\ ,
\end{equation}
with $F$ continuous in the first variable and globally analytic in the
second variable. This gives rise to the perturbed Benjamin-Ono
equation
\begin{equation}
  \label{pertBO}
\partial_tu={\rm H}\partial_x^2 u-\partial_x(u^2)+\epsilon \Pi_0
f(x,\partial_x^{-1}u)\ ,
\end{equation}
where $\Pi_0$ is the projector on states of zero mean:
$$
\Pi_0u:=u-\frac{1}{2\pi}\int_0^{2\pi}u(x)dx\ ,
$$
and $f(x,y):=\partial_yF(x,y)$.

A second example is the Hamiltonian variant of  the damping used in \cite{Gas22} and
is given by 
$$ P(u):=\frac{1}{2}\left(\int_{0}^{2\pi}
u(x)\cos xdx\right)^2+\frac{1}{2}\left(\int_{0}^{2\pi} u(x)\sin xdx\right)^2\ ,
$$
which gives rise to the perturbed BO equation
\begin{equation}
  \label{gassot}
\partial_tu={\rm H}\partial_x^2 u-\partial_x(u^2)+\epsilon\left(\langle
u,\sin(.)\rangle \cos x-\langle
u,\cos(.)\rangle \sin x\right)\ .
\end{equation}
The fact that this perturbation is Hamiltonian is what guarantees the
exponentially long times of stability of the actions.

\section{Proof of Theorem \ref{main1}.}\label{proof}

\subsection{Scheme of the proof}\label{scheme}
We recall that, since we are studying large amplitude solutions, the
frequencies and the resonance relations that they fulfill depend on
the initial datum, so one has to  develop also the so called
geometric part of the proof of Nekhoroshev's theorem.  Here we use
Lochak's approach to the geometric part.

We recall that Lochak's proof is based on the idea of first proving
long time stability of resonant tori and then showing that any initial
datum falls in the domain of stability of some resonant torus, thus
deducing long time stability of every torus.

We now describe more in detail the two steps and illustrate the way
they are developed here for a perturbation of the BO equation.

First of all we fix a resonant finite gap torus, with gaps
$\gamma_1^*,...,\gamma_N^*$, and expand the Hamiltonian of the BO
equation close to it. This is done explicitely in Eq. \eqref{espando}.
One has
$$
H_{BO}(\gamma^*+I)=H_{BO}(\gamma^*)+\sum_{n=1}^N  (n^2-2y_n^*)I_n+\sum_{n\geq N+1}
(n^2-2y_N^*)I_n -H_4(I)\ ,
$$
with $y_n^*$ a linear expression in $\gamma_n^*$ (given explicitely by
\eqref{formule} and \eqref{formule2}) and $H_4$ given by
\eqref{BOH1}. Neglecting the term $H_4$ one gets the expression of the
linearized Hamiltonian at the torus, in particular one gets that the
frequency of motion of such a linearized Hamiltonian are given by 
$$
\left\{\omega_n(\gamma):=n^2-2 y_n(\gamma^*)\right\}_{n=1}^N\ ,\quad
\left\{\omega_n(\gamma):=n^2-2 y_N(\gamma^*)\right\}_{n\geq N+1}\ . 
$$
If one chooses the gaps $\gamma^*$ in such a way that
$y_n(\gamma_*)=\frac{k_n}{q}$ with $k_1,...,k_N,q\in\Z$ then the
linearized dynamics turn out to be periodic with period $q$.

The first important remark by Lochak is that it is particularly simple
and effective to make averages when the unperturbed dynamics is
periodic. So the next step consists in averaging the complete non
integrable system close to a fixed resonant torus. However in order to be able
to perform the subsequent steps of the proof one has to be very
quantitative and to keep into account the size of the neighbourhood of
the torus in which the normal form is valid and the dependence of all
the constants on the period of the unperturbed dynamics. This is done
in Theorem \ref{mainNF}, which in turn is obtained by applying the
main normal form theorem of \cite{Bam99} which is recalled in the
appendix (see Theorem \ref{mainBG}). Theorem \ref{mainNF} requires as
input a precise estimate of the size of the different parts of the
Hamiltonian and of their vector fields in a complex neighbourhood of the
resonant torus. This is obtained in Lemma \ref{Ezero} and Lemma
\ref{omegaf}.

Then, in order to describe the result of the normal form theorem and
the way the normal form is used in order to prove stability of the
resonant tori we have to introduce some notations. We denote by
$h_\omega$ the Hamiltonian generating the linearized dynamics at the
torus:
$$
h_\omega (I):=\sum_{n=1}^N  (n^2-2y_n^*)I_n+\sum_{n\geq N+1}
(n^2-2y_N^*)I_n\ ,
$$
and remark that, since $I_n$ are the difference of the actions and the
gaps of the
resonant torus $\gamma _n^*$, they are not necessarily
positive, so that $h_\omega$ is not positive definite.
The normal form Theorem \ref{mainNF} ensures the existence of a canonical
transformation, defined in a neighbourhood of size $R$ of the resonant
torus, conjugating the Hamiltonian of the perturbed BO equation to a
Hamiltonian of the form
$$
H':=h_\omega-H_4 +Z+\cR\ 
$$
(see \eqref{4.4.1}), where $Z$ has the property that
$\left\{h_\omega,Z\right\}\ =0$ and has the same size as the original
perturbation, namely $\epsilon$, while $\cR$ is a remainder which is
exponentially small in the inverse of 
$$
\mu:=\left (\frac{\epsilon}{R^2}+R^2\right )q\ .
$$ 
Remark that the small parameter $\mu$ contains the size $\epsilon $ of the
perturbation, the size $R$ of the neighbourhood and the period $q$ of the
linearized dynamics at the torus.

We describe now how to use such a normal form to show that the
resonant torus is stable over long times. To this end remark that if
in $H'$ one neglects the exponentially small remainder $\cR$, then one
has two integrals of motion, namely $H' $ and $h_\omega$. If one considers
also the remainder, then $h_{\omega}$ is no more an integral of
motion, but it moves by a small quantity over very long times (while
$H'$ remains an integral of motion). We aim now at exploiting the
convexity of $H_4$ in order to control the motion of each one of the
actions. To this end write the conservation of
energy for $H'$: evaluating at $t=0$ or at a general instant of time,
one has
$$
h_\omega(t)-H_4(t) +Z(t)+\cR(t)=h_\omega(0)-H_4(0) +Z(0)+\cR(0)\ ;
$$
reorganizing the terms and using the triangular inequality, one
gets 
$$
H_4(t)\leq
|h_{\omega}(t)-h_{\omega}(0)|+|Z(t)|+|Z(0)|+|\cR(t)|+|\cR(0)|+H_4(0)\ .
$$ 
Now, all the terms in the right hand side  are small over exponentially long
times, at least if $H_4(0)$ is small, i.e., if the initial datum is
close to the resonant torus. Then the explicit form of $H_4$ allows
us to ensure that all the quantities $s_l$ are small at time $t$, and
therefore that the solution is close to the resonant torus at time
$t$. However the situation is slightly subtler since one needs
closeness in some $H^s$ topology, namely a topology in which one
can ensure the perturbation to be small. However the topology in which
$H_4$ is convex is weaker than an $H^s$ topology. So, in
order to conclude the proof one also uses the almost conservation of
$h_\omega$: the fact that $h_\omega(t)$ is small provides some
additional information and allows us control the distance of the solution from the
resonant torus in the energy norm. This is obtained in
Subsection \ref{stabres}.

After establishing the long time stability of resonant tori one wants
to extend stability to general finite gap tori.

Let us fix a torus whose stability has to be studied. This is
done by fixing
its gaps, say $\gamma_1^0,...,\gamma_N^0$, correspondingly one gets a
sequence of frequency modulations $y_1^0,...,y_N^0$, and one looks for
a resonant sequence $y_1^*,...,y_N^*$ close to it. This is done 
by using Dirichlet theorem, according to which, for any $Q>1$ there exist
$k_1,...,k_N\in\Z$ and $\Z_+\ni q\leq Q$ such that 
$$
\left|y_j^0-\frac{k_j}{q}\right|\leq \frac{1}{qQ^{1/N}}\ ,\quad
\forall j=1,...,N\ . 
$$
By the fact that the map $\gamma\to y$ is invertible, one then
identifies a torus with gaps $\gamma^*$, at a distance
$R\simeq \cO\left(\frac{1}{qQ^{1/N}}  \right)^{1/2}$ from the original one,
on which the linearized dynamics is periodic with period $q$. Thus,
the actions remain close to the actions of the resonant torus for a time
exponentially long with the inverse of 
$$
\left (\frac {\epsilon}{R^2}+R^2\right )q\simeq \frac{\epsilon q}{R^2}+\frac{1}{Q^{1/N}}\lesssim \epsilon Q^{2+\frac 1N}+\frac{1}{Q^{1/N}}\ ,
$$
therefore they also remain close to their initial value. 
Choosing $Q$ in such a way that the two terms of the small parameter
are equal, one gets the wanted exponential estimate valid for all
initial data. 

We emphasize that a few technical points have also to be taken into
account: they are related to the fact that the above conservation
arguments are presented in the coordinates introduced by the
normalizing change of coordinate, which is different for different 
resonant tori. The other point is that one has to
show that the domain in which the normal form is performed is not left
by the initial data one is studying. This is the heart of the proof of
Theorem \ref{stab.res}.

\subsection{Normal form close to resonant tori}\label{resonant}

From now on, in order to take advantage of the analyticity assumption of the Hamiltonian function,
 we work in the complexification of the phase space, which
means that the variables $\xi_n$ have to be considered as independent
of $\bar \xi_n$, so we will consider as phase space $h^1:=h^1_+\oplus h^1_+$
and a point of $h^1$ will be denoted $(\xi,\eta)\equiv ((\xi_n)_{n\geq
  1},(\eta_n)_{n\geq1})$. The real subspace then corresponds to
$\eta_n=\bar\xi_n$.  

\subsubsection{Preliminary estimates}
Following \cite{GK21} we will denote by $\gamma_n:=\xi_n\eta_n$ the
gaps. We now recall some formulae from \cite{GKT22}. In terms of the
Birkhoff coordinates the BO equation takes the form 
\begin{align}
  \label{BOBirkhof}
\dot\xi_n=\im\omega_n(\gamma) \xi_n\ ,\quad \dot
\eta_n=-\im\omega_n(\gamma)\eta_n \ ,\quad n\geq 1
\end{align}
with
\begin{align}
  \label{formule}
\omega_n(\gamma):=n^2-2 y_n(\gamma)\ ,\quad y_n(\gamma):= \sum_{l=1}^ns_l(\gamma)\ ,
\end{align}
and
\begin{equation}
  \label{formule2}
s_l(\gamma):=\sum_{k\geq l}\gamma_k\ , \quad l\geq 1\ .
\end{equation}
By these formulae it is easy to see that the nonlinear correction to
the action to frequency map, namely the map $\gamma\mapsto y$,  is
invertible and its inverse is
\begin{equation}
  \label{ytog}
\gamma_n(y):=2y_{n}-y_{n+1}-y_{n-1} ,\quad n\geq1\ ,
\end{equation}
with the convention $y_0:=0$. We remark that $\omega_n(\gamma)$ are the frequency of motion on the
torus corresponding to the state with actions $\gamma$. 

In particular it follows that for $N$-gap states one has $y_n=y_N$
$\forall n\geq N+1$, and viceversa, given any $N$ dimensional vector
$(y_1,...,y_N)$ with the open condition $y_N>y_{N-1}$ and $2y_n> y_{n+1}+y_{n-1}$, $n=1,\dots, N-1$, one gets a corresponding $N$ gap state.

Fix now a reference $N$ dimensional invariant torus with gaps
$\gamma^*\equiv (\gamma_1^*,...,\gamma_N^*)$ that we will also assume
to fulfill
\begin{equation}
  \label{nozero}
E_m<\gamma_n^*<E_M\ ,\quad \forall n=1,...,N\ ;
\end{equation}
correspondingly we consider the 
``frequency vector'' $y^*\equiv (y^*_1,...,y^*_N) $ defined by
\eqref{formule}, \eqref{formule2}.  
In particular we thus have
\begin{equation}
  \label{y*}
\left|y_n^*\right|\leq{\ \ n\left (N-\frac{n-1}2\right ) E_M }\leq \frac{N(N+1)}{2}E_M\leq N^2E_M\ . 
  \end{equation}
We expand the BO
Hamiltonian at $\gamma^*$. Writing $I_n:=\gamma_n-\gamma^*_n, n=1,\dots, N, $  and
$I_n:=\gamma_n, n\geq N+1$, one gets, having
 set $\gamma_n^*:=0$ for $n\geq N+1$,
 \begin{align}\nonumber
H_{BO}(I+\gamma^*)&=\sum_{n=1}^\infty n^2(\gamma_n^*+I_n)-\sum_{n=1}^\infty [s_n(\gamma^*)+s_n(I)]^2\\ \nonumber
&=\sum_{n=1}^N n^2\gamma_n^*-\sum_{n=1}^N [s_n(\gamma^*)]^2+\sum_{n=1}^\infty n^2I_n -2\sum_{n=1}^\infty s_n(\gamma^*)s_n(I)
- \sum_{n=1}^\infty [s_n(I)]^2\\ \nonumber
&= H_{BO}(\gamma^*)+\sum_{n=1}^\infty  n^2I_n -2\sum_{n=1}^\infty
s_n(\gamma^*)\left (\sum_{k\geq n} I_k\right )-H_4(I)\\ \nonumber
&=  H_{BO}(\gamma^*)+\sum_{n=1}^\infty  n^2I_n -2\sum_{k=1}^\infty I_k\left (\sum_{n\leq k}s_n(\gamma^*)\right )-H_4(I)\\
   \label{espando}
&=  H_{BO}(\gamma^*)+\sum_{n=1}^N  (n^2-2y_n^*)I_n+\sum_{n\geq N+1} (n^2-2y_N^*)\gamma_n -H_4(I)\ .
 \end{align}

We now introduce action angle variables for the first $N$ modes,
namely we introduce variables $(I_n,\phi_n)_{n=1}^N$ by
\begin{equation}
  \label{actionangle}
\xi_n=\sqrt{I_n}{\rm e}^{\im \phi_n}\ ,\quad \eta_n=\sqrt{I_n}{\rm e}^{-\im
  \phi_n}\ ,\quad n=1,..., N\ .
\end{equation}
The variables in the phase space will
now be
$$
((I_n)_{n=1}^N,(\phi_n)_{n=1}^N,(\xi_n)_{n\geq N+1},(\eta_n)_{n\geq N+1})\ .
$$
and the BO-Hamiltonian turns out to be \eqref{espando}
with $$I_n\equiv\gamma_n\equiv \xi_n\eta_n\ ,\quad \forall n\geq N+1\ ,$$
 namely the gaps of index larger than $N$ must be considered as
 functions of the Birkhoff coordinates.

We will call $\cC$ the map that to the variables
$(I,\phi,\left(\xi_n\right)_{n\geq N+1},\left(\eta_n\right)_{n\geq
  N+1} )$ associates the Birkhoff variables, namely
\begin{equation}
  \label{defc}
\cC\left(I,\phi,\left(\xi_n\right)_{n\geq N+1},\left(\eta_n\right)_{n\geq
  N+1} \right):=\left(\left(\xi_n^B\right)_{n\geq 1},\left(\eta_n^B\right)_{n\geq
  1} \right)
\end{equation}
with
\begin{equation}
  \label{defc.1}
\xi^B_n:=\left\{
\begin{matrix}
  \sqrt{I_n+\gamma^*_n}{\rm e}^{\im\phi_n}\ \text{if}\ n=1,...,N\ ,
  \\
  \xi_n \ \text{if}\ n\geq N+1
\end{matrix}
\right.
\\
\eta^B_n:=\left\{
\begin{matrix}
  \sqrt{I_n+\gamma^*_n}{\rm e}^{-\im\phi_n}\ \text{if}\ n=1,...,N\ ,
  \\
  \eta_n \ \text{if}\ n\geq N+1
\end{matrix}
\right.
\end{equation}

Since we work in the complex extension of the phase space the
variables $I,\phi$ will be assumed to be complex.

\begin{definition}\label{def:real}
 A state will be said
to be real if one has
\begin{align}
  \label{real.1}
  I_n\in\R\ ,\quad \phi_n\in\T\,\quad n=1,...,N\ ,
  \\
  \label{real.2}
\eta_n=\bar\xi_n\ ,\quad \forall n\geq N+1\ .
  \end{align}
  \end{definition}
 We now define the norm in the phase
space. It will depend on a parameter $R$ which will be eventually linked to
$\epsilon$. So we define
\begin{equation}
  \label{norma}
\left\|(I,\phi,\xi,\eta)\right\|:=
\sum_{n=1}^{N}\frac{n^2|I_n|}{R}+\sup_{n=1,...,N}R|\phi_n|
+\sqrt{\sum_{n\geq N+1}n^2\left(|\xi_n|^2+|\eta_n|^2\right)}\ .
\end{equation}
{\it From now on the balls will be intended with respect to this
norm.} The complex ball with center $u$ and radius $R_1$ in this norm will be
denoted by $B_{R_1}(u)$. 

A function which will play a fundamental role in the rest of the paper is 
 \begin{equation}
   \label{homega}
h_\omega:=\sum_{n=1}^N(n^2-2y^*_n)I_n+\sum_{n\geq
N+1}(n^2{-2y_N^* })\gamma_n\ , \quad \gamma_n:=\xi_n\eta_n\ .
 \end{equation}
We  now define the domain $\cG$ which will then be extended to the
complex domain in order to apply Theorem \ref{mainBG}.

We fix a parameter $\epsilon_1$ and define
\begin{align}
  \label{G}
\cG&:=\left\{u=(I,\phi,\xi,\eta)\ :\ u\ \text{is \ real}\ ,\ H_4(I)\leq
\epsilon_1^2\ ,\quad h_\omega(I)\leq \epsilon_1\right\}\ ,
\\
\label{complex g}
\cG_R&:=\bigcup_{u\in\cG}B_R(u)\ .
\end{align}
Most of the time, the balls will be taken of radius $R$ with $R$ the
same parameter in the definition of the norm \eqref{norma}.

From now on we will use the notation $a\sleq b$ in order to mean ``there
exists a constant $C$ independent of $\epsilon,\epsilon_1,R$ such that
$a\leq C b$\ ''. 

A first property of the states in $\cG$ is given by the next lemma

\begin{lemma}
  \label{DominioReale}
For $u\equiv(I,\phi,\xi,\eta)\in\cG$ one has
  \begin{align}
    \label{L.0}
    |s_n|\sleq
    \frac{\epsilon_1}{n^2}
    \\
    \label{L.1}
    |I_n|\leq 2\epsilon_1\ ,\quad n=1,...N
    \\
    \label{L.2}
    \sum_{n\geq N+1}n^2\xi_n\eta_n\sleq \epsilon_1\ ,
    \\
       \label{L.3}
    |y_n|\sleq
\epsilon_1   \ .
    \end{align}
\end{lemma}
\proof By the reality of the state one gets $s_n\in\R$, thus
\eqref{BOH1} implies
\begin{align*}
  |s_n|\leq \epsilon_1\ ,\quad\forall n\geq 1
  \\
  |I_n|=|s_n-s_{n+1}|\leq 2\epsilon_1\ ,\ n=1,\dots, N\ ,
  \\
  0\leq \gamma_n\leq \epsilon_1\ ,\ n\geq N+1\ .
  \end{align*}
We come to \eqref{L.2}\eqref{L.3}. First we fix $L$
  such that $(L+1)^2-2y^*_N\geq (L+1)^2/2$. Notice that the size of $L$ is controlled by $E_M$.   Then we have
  \begin{align*}
\frac{1}{2}\sum_{n\geq L+1}n^2\gamma_n\leq \sum_{n\geq
  L+1}(n^2-2y_N^*)\gamma_n=h_\omega -\sum_{n=1}^N (n^2-2y_n^*)I_n-\sum_{n=N+1}^L(n^2-2y_N^*)\gamma_n
\\
\leq
\epsilon_1+\sum_{n=1}^N(n^2+2y_n^*)|I_n|+\sum_{n=N+1}^L(n^2+2y_N^*)\gamma_n
\sleq \epsilon_1\ ,
  \end{align*}
  where we have used $h_\omega \leq \epsilon _1$ on $\mathcal G$.
This implies also $\sum_{n\geq N+1}n^2\gamma_n\sleq \epsilon_1$.

Now, one has, for $l\geq N+1$, 
\begin{equation}
  \label{L.4.1}
s_l=\sum_{k\geq l}\gamma_k=\sum_{k\geq l}\frac{k^2}{k^2}\gamma_k\leq
\frac{1}{l^2}\sum_{k\geq N+1}k^2\gamma_k\sleq \frac{\epsilon_1}{l^2} \ ,
\end{equation}
which clearly holds also for $n\leq N$.
This gives 
\begin{align}
  \label{L.5}
|y_n|\sleq \sum_{l=1}^{n}|s_l|\sleq
\epsilon_1\sum_{l\geq 1}\frac{1}{l^2}\sleq \epsilon_1 \ .
  \end{align}
\qed

\begin{lemma}
  \label{DominioEsteso}
 For $u\in\cG_R$ one has
  \begin{align}
    \label{L.11}
    |I_n|\sleq \epsilon_1+R^2\ ,\quad n=1,...N
    \\
    \label{sn}
|s_n|\sleq    \frac{\epsilon_1+R^2}{n^2}
    \\
    \label{L.12}
    \sum_{n\geq N+1}n^2\gamma_n\sleq\epsilon_1+R^2
    \ ,
    \\
    \label{L.13}
    |y_n|\sleq \epsilon_1+R^2
    \ .
    \end{align}
  \end{lemma}
\proof Recall that 
$$\left\|(I,\phi,\xi,\eta)\right\|:=
\sum_{n=1}^{N}\frac{n^2|I_n|}{R}+\sup_{n=1,...,N}R|\phi_n|
+\sqrt{\sum_{n\geq N+1}n^2\left(|\xi_n|^2+|\eta_n|^2\right)}\ .
$$
First one clearly has $\left|I_n\right|\leq
2\epsilon_1+R^2$. Then, let $u\equiv(I,\phi,\xi,\eta)\in\cG$,
$u'\equiv(I',\phi',\xi',\eta')\in B_R(0)$. By the triangle inequality, we have

\begin{align*}
  \sqrt  {\sum_{n\geq N+1}n^2(|\xi_n+\xi'_n|^2+|\eta_n+\eta'_n|^2) } &\leq
    \sqrt  {\sum_{n\geq N+1}n^2 (|\xi_n|^2+|\eta_n|^2)}+R  \\
  &\sleq \sqrt {\epsilon_1}+R\ . 
\end{align*}
Inserting in \eqref{L.4.1} and in \eqref{L.5} (mutatis mutandis), one
gets the estimate of $s_n$ and of $y_n$.  \qed

\begin{lemma}
  \label{Ezero}
  For $u\in\cG_R$ one has
  \begin{align}
       \label{Ezero.1}
|H_4(u)|\sleq \epsilon_1^2+R^4\ ,
\\
\label{omega0}
\frac{1}{R}\left\|X_{H_4}(u)\right\|\sleq \epsilon_1+R^2\ .
  \end{align}
  \end{lemma}
  \proof The computation of the supremum of $H_4$ is trivial in view
  of  \eqref{sn}. To estimate the vector field of $H_4$, just remark
  that it is given by
  \begin{equation}
    \label{comp}
(0,2y_n,\im 2y_n\xi_n,-\im 2y_n\eta_n) \ ,
    \end{equation}
  with $y_n$ the function of $I_n$ and $\gamma_n$ given by
  \eqref{formule} and \eqref{formule2}. 
  Using the definition of the norm one estimates the norm of
  \eqref{comp} by twice the following quantity,
  $$
R\sup_{n=1,...,N}|y_n|+\sqrt{\left(\sup_{n\geq
    N+1}|y_n|\right)^2\sum_{n\geq
    N+1}n^2\left(|\xi_n|^2+|\eta_n|^2\right) }\ ,
  $$
which, by Lemma \ref{DominioEsteso}, immediately gives \eqref{omega0}. \qed

\begin{lemma}
  \label{omegaf}
 Under the assumptions of Theorem \ref{main1},  provided $R$ and $\epsilon_1$ are small enough, one has
  \begin{align}
    \label{sti.P.1}
    \sup_{u\in\cG_R}\left|P(u)\right|\sleq 1 
    \\
        \label{sti.P.2}
    \sup_{u\in\cG_R}\left\|X_P(u)\right\|\sleq \frac{1}{R}\ .
    \end{align}
  \end{lemma}
\proof First we define the $N$- gap manifold in action-angle-Birkhoff
coordinates by
$$
\cG^B_N:=\left\{ (\xi,\eta)\in h^{1}\ :\ \xi=\bar
\eta\ \text{and}\ \xi_n=\eta_n=0\ \forall n\geq N+1\right\}\ ,
$$
and extend it to the complex domain by defining
$\cG^B_{N,\rho}:=\bigcup_{u\in\cG^B_N}B^{h^{1}}_\rho(u)$.

Then we remark that, by Lemma \ref{DominioReale} one has that
$$\cC(\cG\cap\left\{(\xi_n,\eta_n)=0\ ,\forall n\geq N+1\right\})$$ is a
bounded subset of  $\cG^B_N$, thus, by Lemma
\ref{DominioEsteso}, provided $\epsilon_1$ and $R$ are small enough,
one has that
\begin{equation}
  \label{da c}
\cC(\cG_R)\subset \cG^B_{N,C(R+\sqrt{\epsilon_1})}\ , 
\end{equation}
and furthermore the set in the left hand side  is bounded. Furthermore, if $R$ and
$\epsilon_1$ are small enough, then $\cC(\cG_R)$ is in the domain of
analyticity of $\Phi$ and furthermore $\Phi(\cC(\cG_R))$ is contained
in the domain of analyticity and boundedness of $P$ and of its vector
field. It follows that 
\begin{equation}
  \label{sti.sup}
\sup_{u\in \cG_R}|P(\Phi(\cC(u)))|\sleq 1\ .
\end{equation}
since the supremum is taken over a domain smaller than the domain of
boundedness of $P$. Furthermore, since
    $$
X_{P\circ\Phi}(u)=(d\Phi(u))^{-1}X_P(\Phi(u))\ ,
    $$
one also has
\begin{equation}
  \label{sti.pk}
\sup_{u\in\cG_R}\left\|X_{P\circ\Phi}(\cC(u))\right\|\sleq 1\ .
  \end{equation}

We still have to introduce the action angle variables. To this end we
use $X_{P\circ\Phi\circ\cC}(u)=[d\cC(u)]^{-1}X_{P\circ
  \Phi}(\cC(u))$. We compute now $d\cC$. To
this end remark first that this linear map is the identity on the modes with index
larger than $N$, so consider the smaller indexes: denoting
$(\xi',\eta')^T=d\cC(u)(h_{I},h_\phi)^T$, one has, for $n\leq N$,
\begin{align*}
\xi'_n=\frac{{\rm e}^{\im
    \phi_n}}{2\sqrt{I_n+\gamma_n^*}}h_{I_n}+\im\sqrt{I_n+\gamma_n^*}
{\rm e}^{\im\phi_n} h_{\phi_n}
\\
\eta'_n=\frac{{\rm e}^{-\im
    \phi_n}}{2\sqrt{I_n+\gamma_n^*}}h_{I_n}-\im\sqrt{I_n+\gamma_n^*}
{\rm e}^{-\im\phi_n} h_{\phi_n}
\end{align*}
whose inverse linear map is
\begin{align*}
  h_{I_n}=\left(\xi'_n{\rm e}^{-\im\phi_n}+\eta'_n{\rm e}^{\im\phi_n}\right)
  \sqrt{I_n+\gamma_n^*} 
  \\
  h_{\phi_n}=\left(\xi'_n{\rm e}^{-\im\phi_n}-\eta'_n{\rm e}^{\im\phi_n}\right)
  \frac{1}{\im \sqrt{I_n+\gamma_n^*} }\ .
\end{align*}
To evaluate the norm of this map, we compute
$$
\left\|
    [d\cC(u)]^{-1}(\xi',\eta')^T\right\|=\frac{1}{R}\sum_{n=1}^{N}n^2\left|h_{I_n}\right|+R\sup_{n=1,...,N}\left|h_{\phi_n}\right|\ .
    $$
    If $|I_n|\leq \gamma_n^*/2$, which is ensured by $R^2<E_m/8$, the
    second term in the right hand side is bounded uniformly in $R$. The first term is bounded
    by a constant times $R^{-1}$, so one gets \eqref{sti.P.2}.
\qed

\subsubsection{Normal form}

From now on we take
\begin{equation}
  \label{epsilon1}
\epsilon_1=R^2\ ,
\end{equation}
so that, on $\cG_R$, one has
\begin{align}
  \label{9.1}
  \left|H_4\right|\sleq R^4\ ,\quad \left\|X_{H_4}\right\|\sleq R^3
  \\
  \label{9.2}
  \left|\epsilon P\right|\sleq\epsilon\ ,\quad \left\|X_{\epsilon
    P}\right\|\sleq \frac{\epsilon}{R}\ . 
\end{align}

We now fix a resonant torus, make a normal form close to it and use
$h_{\omega}$ and the energy as Lyapunov functions in order to prove
long time stability of this torus. For definiteness we avoid here the
use of the symbol $\sleq$.

\begin{theorem}
  \label{mainNF} Under the assumptions of Theorem \ref{main1}, assume that $$
  y^*_n=\frac{k_n}{q}\ ,\quad n=1,...,N\ ,
  $$
with $k_1,...,k_n,q \in\Z_+$, 
  then there exist positive constants $\mu_*,K_1,...,K_5$ (independento of
$k_n$ and $q$) such that, if
  \begin{equation}
    \label{NF.1}
\mu:=(R^2+\frac{\epsilon}{R^2})q<\mu_*/2\ ,
  \end{equation}
then there exists an analytic canonical transformation
$\cT:\cG_{\erremu}\to\cG_{R}$ with $\cT(\cG_{\erremu})\supset \cG_{\erremuu}$  which is close to identity, namely
\begin{align}
  \label{close.1}
\sup_{u\in\cG_{\erremu}}\left\|u-\cT(u)\right\|\leq K_1\frac{\epsilon}{R}q\ ,
\\
  \label{close.2}
\sup_{u\in\cG_{\erremuu}}\left\|u-\cT^{-1}(u)\right\|\leq K_1\frac{\epsilon}{R}q\ ,
\end{align}

\begin{equation}
\end{equation}
and is such that
\begin{equation}
  \label{4.4.1}
H\circ\cT=h_\omega-H_4 +Z+\cR\ ,
\end{equation}
with
\begin{itemize}
\item $Z$ is analytic on $\cG_{\erremu}$ toghether with its Hamiltonian
  vector field and fulfills
  \begin{equation}
    \label{NF.4}
\sup_{u\in\cG_{\erremu}}\left|Z(u)\right|\leq  K_2\epsilon\ ,\quad
\sup_{u\in\cG_{\erremu}}\left\| X_Z(u)\right\|\leq K_3\frac{\epsilon}{R}\ ,  
  \end{equation}
  furthermore it is in normal form, namely one has
  \begin{equation}
    \label{nom}
\left\{h_{\omega},Z\right\}\equiv dh_{\omega}(u)X_{Z}(u)=0\ ,\forall
u\in \cG_{\erremu}\ . 
  \end{equation}
\item $\cR$ is analytic on $\cG_{\erremu}$ toghether with its Hamiltonian
  vector field and fulfills
  \begin{align}
      \label{stimeZ1.1}
\sup_{u\in\cG_{\erremu}}\left|\cR(u)\right|\leq
K_4\epsilon\exp\left(-\frac{\mu_*}{\mu}\right)\ ,\\
  \label{st.ff}
\sup_{u\in\cG_{\erremu}}\left\| X_{\cR}(u)\right\|\leq K_5
\frac{\epsilon}{R} \exp\left(-\frac{\mu_*}{\mu}\right) \ , 
  \end{align}
\end{itemize}
\end{theorem}
\proof The theorem is a direct application of Theorem \ref{mainBG} of
the appendix: it is obtained by inserting the values of the constants
obtained from Lemmas \ref{Ezero}, \ref{omegaf}. We just remark
that we have
$$
E^{\sharp}\simeq
\frac{\epsilon}{R^2}\frac{3\epsilon+R^4}{9\frac{\epsilon}{R^2}
  +5R^2}=\epsilon\frac{ 3\epsilon+R^4}{9\epsilon+5R^4}\simeq
\epsilon\ .
$$\qed

Applying Remark \ref{deltaF}, one also gets the
following result.
\begin{remark}
  \label{deltahomega}
    By the above preliminary estimates one has
  \begin{equation}
    \label{homega1}
\sup_{u\in\cG_R}\left|h_\omega(u)\right|\sleq R^2\ ,
  \end{equation}
  so that, by \eqref{deltaF.1}, we have
  \begin{align}
  \label{deltaomega}
    \sup_{u\in\cG_{\erremu}}\left|h_\omega(u)-h_\omega(\cT(u))\right|\sleq
 \epsilon q
    \\
    \label{3.40}
    \sup_{u\in\cG_{\erremu}}\left|H_4 (u)-H_4(\cT(u))\right|\sleq R^2 \epsilon q\ ,
  \end{align}
  and similarly with $\cT$ replaced by $\cT^{-1}$.
\end{remark}

\subsection{Stability of resonant tori}\label{stabres}

From the above theorem and remark one gets the stability of resonant tori.
\begin{theorem}
  \label{stab.res}
There exist positive constants $\mu_*\ll1$ and $K,\tilde K\gg 1$ with the
following properties: assume
\begin{equation}
  \label{varie}
\mu<\mu_*/2\ ,\quad \epsilon q\leq\frac{R^2}{K}\ ,\quad \epsilon\leq
\frac{R^4}{\tilde K}\ ,
\end{equation}
and consider initial data $u_0$ with
\begin{equation}
  \label{ini.data}
h_\omega(u_0)\leq \frac{R^2}{2}\ ,\quad H_4 (u_0)\leq
\left(\frac{R^2}{2}\right)^2\ .
\end{equation}
Then, for all times $t$ fulfilling
\begin{equation}
  \label{tempi}
|t|\leq\frac{R^4}{K\epsilon}\exp\left(\frac{\mu_*}{\mu}\right)\ ,
\end{equation}
the solution of the perturbed BO equation \eqref{Hamilton}, \eqref{hperturbata}, fulfills
\begin{equation}
  \label{sti.tempo}
h_\omega(u(t))\leq R^2\ ,\quad H_4 (u(t))\leq R^4\ ,
\end{equation}
and in particular, in view of \eqref{G} and \eqref{epsilon1},  $u(t)\in\cG$.
\end{theorem}
\proof We proceed by a bootstrap argument. Assume that there exists $t$ satisfying \eqref{tempi} and such that \eqref{sti.tempo} does not hold. Denote by 
$t^*$ a time of minimal absolute value  satisfying \eqref{tempi} and such that  
\begin{equation}\label{tstar}
h_\omega(u(t^*))= R^2\quad {\rm or}\quad H_4 (u(t^*))= R^4\ ,
\end{equation} 
and let us look for a contradiction. \\
First, in view of \eqref{tstar}, we can make the canonical transformation $\cT$ provided by Theorem \ref{mainNF} near the trajectory for $|t|\leq |t^*|$, and we set, for $|t|\leq |t^*|$, 
$$u=\cT(u')\ ,\ h_\omega(t):=h_\omega(u(t))\ ,\ h'_\omega(t):=h_\omega(u'(t))\ ,$$ and similarly for $H_4 $. From \eqref{ini.data} and \eqref{deltaomega} we have,
with a constant $C$ that can change from line to line,
\begin{equation*}
\left|h'_{\omega}(0)\right|\leq
\frac {R^2}{2}+C \epsilon q
\end{equation*}
and therefore, since 
\begin{equation}
  \label{spostoomega}
\left|\dot
h_\omega'\right|=\left|\left\{h_{\omega},\cR\right\}\right|=\left|dh_\omega
X_{\cR}\right|\leq
C\frac{R^2}{R}\frac{\epsilon}{R}\exp\left(- \frac{\mu_*}{\mu}\right)
\end{equation}
one has
\begin{align}
h_\omega(t)&\leq
|h_{\omega}(t)-h'_{\omega}(t)|+|h'_{\omega}(t)-h'_{\omega}(0)|+|h'_{\omega}(0)|
\\
&\leq2\epsilon q C+ C\epsilon\exp\left(-\frac{\mu_*}{\mu}\right)|t|+   \frac{R^2}{2}
\end{align}
which, provided
\begin{equation}\label{sti.homega.1}
2\epsilon q C\leq \frac{R^2}{8}\ ,\quad
C\epsilon\exp\left(-\frac{\mu_*}{\mu}\right)|t|\leq \frac{R^2}{8} 
\end{equation}
is not bigger  than $3R^2/4$. In particular $h_\omega (t^*)\leq 3R^2/4$. Notice that the second estimate of
\eqref{sti.homega.1} is ensured by \eqref{tempi} and the first one by
\eqref{varie}. 

We come to $H_4 $. Exploiting the conservation of the Hamiltonian \eqref{4.4.1}, one
gets
$$
H_4 '(t)-H_4'(0)=h'_{\omega}(t)-h'_{\omega}(0)+Z(t)-Z(0)+\cR(t)-\cR(0)\ ,
$$
which gives
$$
\left|H_4 '(t)-H_4'(0)\right|\leq
C\epsilon\exp\left(-\frac{\mu_*}{\mu}\right)|t| +2C\epsilon + 2
K_4\epsilon \exp\left(-\frac{\mu_*}{\mu}\right)\ ,
$$
taking also into account \eqref{3.40}
one gets
\begin{align*}
H_4(t)\leq H_4(0)+\left|H_4(0)-H_4'(0)\right|+
\left|H_4(t)-H_4'(t)\right|+ \left|H_4'(t)-H_4'(0)\right|
\\
\leq \frac{R^4}{4}+C R^2\epsilon q+C\epsilon+
C\epsilon\exp\left(-\frac{\mu_*}{\mu}\right)|t|\ ;
\end{align*}
thus, if each term (but the first) is smaller than $R^4/16$ one gets
that $H_4(t^*)\leq 7R^4/16$.  Again, this is ensured by \eqref{varie} and
\eqref{tempi}.\\
Summing up, we have proved 
$$ h_\omega (t^*)\leq 3R^2/4\ ,\ H_4(t^*)\leq 7R^4/16\ ,$$
which obviously contradicts \eqref{tstar}. 
\qed

\subsection{Stability of finite gap tori}\label{stabfinite}

\noindent {\it Proof of Theorem \ref{main1}.} We take now an initial
datum $\xi_n^0$ fulfilling \eqref{ini.1}. We also assume, instead of
\eqref{ini.2}
\begin{equation}
  \label{ini.3}
\sum_{n\geq N+1}n^2|\xi_n^0|^2\leq {\epsilon_2}\ ,
  \end{equation}
with $\epsilon_2$ to be determined later. We
look for a resonant torus close to it. Consider
$\gamma_n^0:=|\xi_n^0|^2$, for $n=1,..., N$. We approximate the torus
with gaps $\gamma_n^0$ by a resonant torus. To this end define $y_n^0$
by \eqref{formule} so that $|y_n^0|\leq N^2E_M$ and use Dirichlet's
theorem to approximate $y_n^0$ by a rational vector. Indeed Dirichlet's
theorem ensures that for every $Q\geq1$ there exist integers $k_1,...,k_N,
q\in \Z$ with $1\leq q\leq Q$ such that
\begin{equation}
  \label{diri}
\left|y_n^0-\frac{k_n}{q}\right|\leq \frac{1}{qQ^{1/N}}\ .\quad n=1,...,N
  \end{equation}
We now define $y_n^*:=k_n/q$ and correspondingly $\gamma_n^*$ by
\eqref{y*} and \eqref{ytog}. Choosing $Q\geq (N^2E_M)^{-N}$, which is obviously granted by the choice \eqref{Qbis} below, we have
\begin{equation}
  \label{diri.1}
|y_n^*|\leq 2N^2 E_M\ ,\quad \left|\gamma_n^*-\gamma_n^0\right|\leq \frac{4}{qQ^{1/N}}\ ,\quad
|s_n(\gamma^0)-s_n(\gamma^*)| \leq \frac{2}{qQ^{1/N}}\ .
\end{equation}
We can now compute $h_\omega(0)$:
passing to the action variables $I_n$, we have
\begin{eqnarray*}
\sum_{n=1}^{N}(n^2+2y_n^*)\left|I_n\right|&\leq
& \sum_{n=1}^{N}(n^2+4N^2E_M) \frac{4}{qQ^{1/N}}\\
&\leq &(N^3+4N^3E_M)
\frac{4}{qQ^{1/N}}\leq \frac{20N^3E_M}{qQ^{1/N}}\ ,
\end{eqnarray*}
{assuming $E_M\geq 1$, from which one also has
  $$
\sum_{n\geq N+1}\left|y_N^*\right|\left|\xi_n^0\right|^2\leq
\frac{\left|y_N^*\right|^2}{(N+1)^2}\sum_{n\geq
  N+1}n^2\left|\xi_n^0\right|^2\leq \frac{2N^2E_M}{(N+1)^2}\epsilon_2<2E_M\epsilon_2\ .
  $$
This gives }
\begin{align}
  \label{diri.2}
h_\omega=&\sum_{n=1}^{N}(n^2-2y_n^*)I_n+\sum_{n\geq
  N+1}{(n^2-2y_{N}^*)}\left|\xi_n^0\right|^2
\\
\leq 
&\frac{20N^3E_M}{qQ^{1/N}}+\epsilon_2{+4E_M\epsilon_2}\ .  
\end{align}
Taking
{\begin{equation}
  \label{epsilon2}
\epsilon_2:={\frac{10N^3}{qQ^{1/N}}}\ ,
\end{equation}}
one gets
\begin{equation*}
  \epsilon_2(1+4E_M)\leq \frac{20N^3E_M}{qQ^{1/N}}\ ,
\end{equation*}
and therefore 
$$
h_\omega(0)\leq \frac{40N^3E_M}{qQ^{1/N}}\ .
$$
We come to $H_4$. Working as in \eqref{L.4.1}, one gets
\begin{equation}
  \label{diri.4}
H_4 (0)\leq \sum_{n=1}^{N}\left(\frac{2}{qQ^{1/N}}\right)^2 +K \epsilon_2^2
  \end{equation}
where $K:=\sum_{n\geq N+1}n^{-4}\leq 1$. So, observing that  $\epsilon_2\geq
\sqrt N\frac{2}{qQ^{1/N}} $ is granted by \eqref{epsilon2}, one gets
\begin{equation}
  \label{hat h}
H_4 (0)\leq (K+1)\epsilon_2^2\leq 2\epsilon_2^2\ .
  \end{equation}
In order to ensure \eqref{ini.data}, we take 
\begin{equation}
  \label{epsilon11}
R^2:= \frac{40N^3{E_M}}{qQ^{1/N}}\ .
  \end{equation}
 We now aim at applying Theorem \ref{stab.res}. Then we would like to  optimize the
 value of $Q$ in order to maximize the time of validity of the
 estimates. Inserting \eqref{epsilon11} in \eqref{NF.1} we get (with a
 suitable $C$)
 \begin{equation}
   \label{diri.14}
\mu\leq C\left(\frac{1}{qQ^{1/N}}+\epsilon qQ^{1/N} \right)q\leq C\left(\frac{1}{Q^{1/N}}+\epsilon Q^{2+1/N} \right) \ .
 \end{equation}
 This would lead to choose $Q$ by imposing the two terms in brackets to be
 equal. This would give 
 \begin{equation}
   \label{Q}
Q=\epsilon^{-\frac{N}{2(N+1)}}\ ,
 \end{equation}
but we also have to ensure the validity of the last inequality in \eqref{varie},
while \eqref{Q} can only ensure $R^4\geq (40 N^3 E_M)^2\epsilon$,
which is not necessarily bigger than $\tilde K\epsilon$. For this
reason we take
\begin{equation}
  \label{Qbis}
Q:=\left(\max\left\{\frac{\tilde
  K}{(40N^3E_M)^2}, 1 \right\}\epsilon \right)^{-\frac{N}{2(N+1)}},
  \end{equation}
so that
  \begin{equation}
   R^2=\frac{40N^3{E_M}}{qQ^{1/N}}\simeq
\frac{\epsilon^{\frac{1}{2(N+1)}}}{q}\sleq \epsilon^{\frac{1}{2(N+1)}}\ ,\quad \mu =\left (R^2+\frac{\epsilon}{R^2}\right )q\simeq
\epsilon^{\frac{1}{2(N+1)}}\ .
  \end{equation}
 In particular, concerning the second estimate, we remark that, by the
 Dirichlet theorem, $q$ is smaller than $Q$, but it can be of order 1.\\
  Inserting in the different estimates one concludes the proof. \qed

\appendix
\section{Normal form close to resonant states}\label{Bam99}

In this appendix we recall the normal form theorem from \cite{Bam99} that we use in Theorem \ref{mainNF} in order to
prove stability of resonant finite gap solutions in Theorem \ref{stab.res}.

Let $\cP$ be a weakly symplectic space that we will also assume to be
a Hilbert space. Given a smooth Hamiltonian function $H$ we will
denote by $X_H$ the corresponding Hamiltonian vector field.

First we define the family of domains in which the normal form theorem
holds.  Consider the complexification $\cP^{\C}$ of $\cP$ and let
$\cG\subset \cP$ be a (real) domain. For $R>0$ we consider
the domain
$$
\cG_{R}:=\bigcup_{u\in\cG}B_{R}(u)
$$
where $B_{R}(u)\subset\cP^{\C}$ is the open complex ball of radius $R$
and center $u$.

Then we consider a Hamiltonian of the form
\begin{equation}
  \label{a.1}
H(u)=h_{\omega}(u)+\hat h(u)+P(u)\ ,
\end{equation}
where $h_\omega$ will be considered as the unperturbed Hamiltonian and
will be assumed to generate a periodic flow, $\hat h(u)$ will be assumed to
Poisson commute with $h_\omega$ while $P$ will be the true
perturbation.

Precisely we assume that

\begin{itemize}

\item[H0)] $h_\omega$ is a non homeogeneous smooth quadratic polynomial
on $\cP^{\C}$. 
  
\item[H1)] There exists a dense subset $D\subset\cP^\C$ where the
  vector field $X_{h_\omega}$ is defined. By H0 $X_{h_{\omega}}$ is a
  linear affine operator that we assume to generate a flow $\Phi^t$ of
  unitary operators.

\item[H2)] The flow $\Phi^t$ is periodic of period $2\pi/\omega$, namely
  \begin{equation}
    \label{periodico}
\Phi^{t+\frac{2\pi}{\omega}}=\Phi^t\ .
  \end{equation}
\end{itemize}


Concerning $\hat h$ we assume
\begin{itemize}
\item[Q1)] $\hat h$ extends to a bounded analytic function on
  $\cG_{R}$. We denote by $E_0$ a constant such that
  \begin{equation}
    \label{E0}
\sup_{u\in\cG_{R}}|\hat h(u)|\leq E_0\ .
  \end{equation}

  \item[Q2)]  The Hamiltonian vector field $X_{\hat h}$ of $\hat h$ extends to a
  bounded  analytic function from $\cG_{R}$ to $\cP^{\C}$. We denote
  by $\omega_0$ a constant such that
\begin{equation}
    \label{E0.1}
\frac{1}{R}\sup_{u\in\cG_{R}}{\|X_{\hat h}(u)\|}\leq \omega_0\ .
\end{equation}
\item[Q3)] $\hat h$ Poisson commutes with $h_\omega$, namely one has
  $$
\left\{h_\omega,\hat h\right\}(u):=dh_{\omega}(u)X_{\hat h}(u)=0\ ,\quad \forall
u\in\cG_{R} \ .
  $$
  \end{itemize}
Finally, concerning the perturbation $P$ we assume 

\begin{itemize}
\item[P1)] $P$ extends to a bounded analytic function on
  $\cG_{R}$. We denote by $E$ a constant such that
  \begin{equation}
    \label{E0.2}
\sup_{u\in\cG_{R}}|P(u)|\leq E\ .
  \end{equation}

  \item[P2)]  The Hamiltonian vector field $X_P$ of $P$ extends to a
  bounded  analytic function from $\cG_{R}$ to $\cP^{\C}$. We denote
  by $\omega_P$ a constant such that
\begin{equation}
    \label{E0.3}
\frac{1}{R}\sup_{u\in\cG_{R}}{\|X_P(u)\|}\leq \omega_P\ .
\end{equation}
\end{itemize}

\begin{theorem}
[Theorem 4.4 of \cite{Bam99} with $f_1\equiv r
      \equiv g\equiv 0$]   \label{mainBG} Denote
  \begin{align}
     \label{mu}
\mu:={2^6 {\rm e}\pi}\frac{\omega_0+\omega_P}{\omega}\ ,\quad E^{\sharp}:= 
\max\{E,\frac{2\omega_P(3E+E_0)}{9\omega_P+5\omega_0}\}\ ,
  \end{align}
  then the following holds true: if
  \begin{equation}
    \label{mupiccolo}
\mu<\frac{1}{2}\ ,
  \end{equation}
then there exists an analytic canonical transformation
$\cT:\cG_{\erremu}\to\cG_{R}$ with $\cT(\cG_{\erremu})\supset
\cG_{\erremuu}$, which is close to identity, namely
\begin{equation}
  \label{close}
\sup_{u\in\cG_{\erremu}}\left\|u-\cT(u)\right\|\leq 4\pi
R\frac{\omega_P}{\omega}\ ,\quad \sup_{u\in\cG_{\erremuu}}\left\|u-\cT^{-1}(u)\right\|\leq 4\pi
R\frac{\omega_P}{\omega}
\end{equation}
and is such that
\begin{equation}
  \label{4.4}
H\circ\cT=h_\omega+\hat h+Z+\cR\ ,
\end{equation}
with
\begin{itemize}
\item $Z$ is analytic on $\cG_{\erremu}$ together with its Hamiltonian
  vector field and fulfills
  \begin{equation}
    \label{stimeZ}
\sup_{u\in\cG_{\erremu}}\left|Z(u)\right|\leq 2E^\sharp\ ,\quad
\frac{1}{R}\sup_{u\in\cG_{\erremu}}\left\| X_Z(u)\right\|\leq 2\omega_P\ , 
  \end{equation}
  furthermore it is in normal form, namely one has
  \begin{equation}
    \label{nom1}
\left\{h_{\omega};Z\right\}(u)\equiv dh_{\omega}(u)X_{Z}(u)=0\ ,\forall
u\in \cG_{\erremu}\ . 
  \end{equation}
\item $\cR$ is analytic on $\cG_{\erremu}$ together with its Hamiltonian
  vector field and fulfills
  \begin{align}
     \label{stimeZ1}
\sup_{u\in\cG_{\erremu}}\left|\cR(u)\right|\leq E^\sharp \exp\left(-\frac{1}{\mu}\right)\ ,
\\
\label{stime33}
{\frac 1R } \sup_{u\in\cG_{\erremu}}\left\| X_{\cR}(u)\right\|\leq  {\rm e}\, \omega_P\exp\left(-\frac{1}{\mu}\right) \ . 
  \end{align}
\end{itemize}
\end{theorem}

\begin{remark}
  \label{deltaF}
By the Cauchy estimate of the differential and \eqref{close}, if $F$ is analytic {from  $\cG_R$ to $\cP^{\C}$}, then one has
\begin{align}
  \label{deltaF.1}
\sup_{u\in\cG_{\erremu}}\left\|F(u)-F(\cT(u))\right\|\leq
8\pi\frac{\omega_P}{\omega} \sup_{u\in\cG_R}\left\| F(u)\right\|\ ,
\\
\label{deltaFu}
\sup_{u\in\cG_{\erremuu}}\left\|F(u)-F(\cT^{-1}(u))\right\|\leq
8\pi\frac{\omega_P}{\omega} \sup_{u\in\cG_R}\left\| F(u)\right\|\ ,
\end{align}
\end{remark}

\bibliography{biblio}

\bibliographystyle{amsalpha}

\end{document}